\documentclass[10pt,leqno,preprint]{imsart}
\usepackage{amssymb}
\usepackage{amsmath}
\numberwithin{equation}{section}

\def\qed{\ \vrule height4pt width4pt depth0pt}

\newcommand{\vanish}[1]{}

\newtheorem{theorem}{Theorem}

\newcommand{\beas}{\begin{eqnarray*}}
\newcommand{\enas}{\end{eqnarray*}} 
\newcommand{\bea}{\begin{eqnarray}}
\newcommand{\ena}{\end{eqnarray}}

\newcommand{\SSS}{{\mathcal S}}

\newcommand{\one}{{\mathbf 1}}

\begin{document}
\begin{frontmatter}
\title{An Exposition of G\"otze's Estimation of the Rate of
Convergence in the Multivariate Central Limit Theorem}
\runtitle{Stein's Method for Multivariate CLT}

\begin{aug}

\author{\fnms{Rabi}\snm{ Bhattacharya}
\thanksref{t1}
\ead[label=e1]
{rabi@math.arizona.edu}}
 
 \and
 \author{ \fnms{Susan} \snm{Holmes}
\thanksref{t2}
 \ead[label=e2]
 {susan@stat.stanford.edu}
 \ead[label=u2,url]{http://www-stat.stanford.edu/\~{}susan}
 }
 
\thankstext{t1}{Supported by NSF DMS-0806011}
\thankstext{t2}{Supported by NSF DMS-02 41246 and NIH R01GM086884-02 }
\runauthor{Bhattacharya and Holmes}

\affiliation{University of Arizona\thanksmark{m1} and Stanford University\thanksmark{m2}}

\address{
Department of Mathematics \\
Math 603\\
University of Arizona \\
\printead{e1}
}

\address{
Department of Statistics\\
Sequoia Hall\\
 Stanford University\\
CA 94305\\
\printead{e2}\\
\printead{u2}}
\end{aug}

\maketitle

\begin{abstract}
We provide an explanation of the main ideas underlying
G\"otze's main result in
\cite{gotze}
using 
Stein's method.
We also provide  detailed derivations of various intermediate estimates.
Curiously, we are led to a different dimensional dependence of the constant than
that given in \cite{gotze}.

We would like to dedicate this to Charles Stein on the occasion of his 90th birthday.
\end{abstract}



\end{frontmatter}

\newpage
\section{Introduction}
In his article G\"otze\cite{gotze} used Stein's method to provide
an ingenious derivation of the Berry-Esseen type bound for the class of Borel convex subsets of
${\mathbb R}^k$ in the context of the classical multivariate central limit theorem (CLT).
This approach has proved fruitful in deriving error bounds for the CLT under certain structures of dependence as well (see Rinott and Rotar \cite{rinottrotar}). Our view and elaboration of
 G\"otze's proof resulted from a collaboration between the authors and were first presented
 in a seminar
 at Stanford given by the first author in the summer of 2000. The authors wish to
 thank Persi Diaconis for pointing out the need for a more
 readable account of G\"otze's  result than given in his original work.
 
 After an explanation of the general method in Section 1, detailed derivations of various estimates are given in Sections 2-4 in terms that would be reasonably familiar to probabilists. Except for the smoothing inequality in Section 4, which is fairly standard, complete proofs are given.
 
 Recently Raic\cite{raic} has followed essentially the same route as G\"otze, but in greater detail, in deriving G\"otze's  bound. It may be pointed out that we were unable to verify the dimensional
 dependence $O(k)$ in \cite{gotze},\cite{raic}.
 Our derivation provides the higher order dependence
 of the error rate on $k$, namely $O(k^{\frac{5}{2}})$.
 This rate can be reduced to
  $O(k^{\frac{3}{2}})$
  using an inequality of Ball \cite{ball}. The best order of dependence known, namely,
$O(k^{\frac{1}{4}})$ is given by Bentkus\cite{bentkus}, using a different method, which would be difficult
to extend to dependent cases.

As a matter of notation, the constants $c$, with or without subscripts are absolute constants. The $k$-dimensional standard Normal distribution is denoted by ${\mathcal N}(0,{\mathbb I}_k)$
as well as $\Phi$, with density $\phi$.  
\subsection{The Generator of the ergodic Markov process as a Stein operator.}
\label{intro}
Suppose $Q$ and $Q_0$ are two probability measures on a measurable space
$(S,\SSS)$ and h is integrable (with regards to $Q$ and $Q_0$). Consider the problem of estimating
			   \begin{equation}
Eh-E_0h\equiv \int h dQ -\int h dQ_0.
\label{eq1.1}
\end{equation}
A basic idea of  Stein\cite{steinbook} 
(developed in some examples in \cite{diaconisholmes}
and \cite{holmes})
is\begin{description}
\item[(i)] to find an invertible map $L$ which maps ``nice'' functions on $S$ into the 
\underline{kernel} or \underline{null space} of $E_0$, 
\item[(ii)]
to find a perturbation of $L$, say $L_\alpha$, which 
 maps ``nice'' functions on $S$ into the
\underline{kernel} or \underline{null space} of $E$,
\item[(iii)]to estimate \ref{eq1.1}
using the identity
\end{description}
			   			   \begin{equation}
		Eh-E_0h=ELg_0=E(Lg_0-L_\alpha g_\alpha)
\label{eq1.2}
\end{equation}
where $$g_0\equiv L^{-1}(h-E_0h),\qquad
g_\alpha\equiv L_\alpha^{-1}(h-Eh).$$

One way to find $L$ is to consider an ergodic Markov process
$\{X_t: t\geq 0\}$ 
on $S$ which has $Q_0$ as it's invariant distribution, and let $L$ be its
\underline{generator}:
\bea
\label{eq:3.1}
Lg=\lim_{t\downarrow 0} \frac{T_tg-g}{t}, \qquad g \in {\mathcal D}_L
\ena
where the limit is in $L^2(S,Q_0)$ , and $$(T_tg)(x)=E\left[g(X_t)|X_0=x\right],$$
or in terms of the transitions probability $p(t;x,dy)$
of the Markov process 
$\{X_t: t\geq 0\}$ ,
\bea
(T_tg)(x)=\int_S g(y) p(t;x,dy)\qquad (x\in S,t>0).
\ena
Also ${\mathcal D}_L$ is the set of $g$ for which the limit in 
(\ref{eq:3.1}) exists.
By the Markov (or, semigroup) property, $T_{t+s}=T_tT_s=T_sT_t$, so that
\bea
\frac{d}{dt}T_tg=\lim_{s\downarrow 0} \frac{T_{t+s}g-T_tg}{s}
=\lim_{s\downarrow 0} \frac{T_t(T_{s}g-g)}{s}=T_tLg.
\ena
Since $T_tT_s=T_sT_t, T_t$ and $L$ commute so that
\bea
\frac{d}{dt}T_tg=LT_tg.
\ena
Note that invariance of $Q_0$ means $ET_tg(X_0)=Eg(X_0)=\int gdQ_0$, if the 
distribution of $X_0$ is $Q_0$. This implies that, for every
$g\in {\mathcal D}_L$, $ELg(X_0)=0$, or
$$
\int_SLg(x)dQ_0(x)=0,
\left[ ELg(X_0)=E(\lim_{t\downarrow 0} \frac{T_tg(X_0)-g(X_0)}{t})=
\lim_{t\downarrow 0} \frac{ET_tg(X_0)-Eg(X_0)}{t}
\right]
$$
That is, $L$ maps ${\mathcal D}_L$ into the set 
$ 1^{\bot}$ 
of mean zero functions in $L^2(S,Q_0)$. It is known that the range of $L$ is dense in 
$ 1^{\bot}$ 
and if $L$ has a spectral gap,
then the range of $L$ is all of $1^{\bot}$.
In the latter case $L^{-1}$ is well defined on
 $ 1^{\bot}$ (kernel of $Q_0$) and is bounded on it
(\cite{bhattacharya}).

Since $T_t$ converges to the identity operator as $t\downarrow 0$
one may also use $T_t$ for small $t>0$ to smooth the target function
$\tilde{h}=h-\int h dQ_0$. For the case of a diffusion
$\{X_t: t\geq 0\}$, 
$L$ is a differential operator and even non smooth functions such
as  $\tilde{h}={\mathbf 1}_B-Q_0(B) (h={\mathbf 1}_B)$
are immediately made smooth by applying $T_t$. One may then use the approximation
to $\tilde{h}$ given by
\bea
T_t\tilde{h}=L(L^{-1} T_t \tilde{h})=L\psi_t, \mbox{ with }\psi_t=L^{-1}T_t\tilde{h},
\ena
and then estimate the error of this approximation by a ``smoothing inequality'', especially
if $T_t\tilde{h}$ may be represented as a perturbation by convolution.
For several perspectives and applications of Stein's method see 
\cite{barbour}, \cite{diaconisholmes},\cite{holmes},\cite{rinottrotar}.

\subsection*{1(b) The Ornstein-Uhlenbeck Process and its Gausssian invariant Distribution}
The Ornstein-Uhlenbeck (O-U)
process is governed
by the \underline{Langevin equation} (see, e.g. \cite[pp. 476, 597, 598]{bhattacharyawaymire})
\bea
dX_t=-X_tdt+\sqrt{2}dB_t
\ena
where
$\{B_t: t\geq 0\}$ is a $k$-dimensional standard Brownian
motion.
Its transition density is
\bea
\hskip0.5cm
p(t;x,y)=\prod_{i=1}^k \left[2\pi (1-e^{-2t})
\right]^{-\frac{1}{2}}\exp\{-\frac{(y_i-e^{-t}x_i)^2}{2(1-e^{-2t})}\}
\qquad
x=(x_1,\ldots,x_k),y=(y_1,\ldots,y_k).
\ena
This is the density of a Gaussian (Normal) 
distribution with mean vector
$e^{-t}x$ and dispersion matrix $(1-e^{-2t}){\mathbb I}_k$
where 
${\mathbb I}_k$ 
is the $k \times k$ identity matrix. 
One can check (e.g., by direct differentiation) that the Kolmogorov backward equation
holds:
\bea
\hskip0.5cm
\frac{\partial p(t;x,y)}{\partial t}
=
\sum_{i=1}^k
\frac{\partial^2 p(t;x,y)}{\partial x_i^2}
-\sum_{i=1}^k
x_i
\frac{\partial 
p(t;x,y)
}{\partial x_i}
=\Delta p-x\bullet \nabla p=Lp, \mbox{ with }L\equiv\Delta-x\bullet \nabla
\ena
where $\Delta$ is the Laplacian and $\nabla=grad$.
Integrating both sides w.r.t. $h(y)dy$
we see that $T_th(x)=\int h(y) p(t;x,y) dy$ satisfies
\bea
\frac{\partial}{\partial t}
T_th(x)=\Delta T_th(x)-x\bullet\nabla T_th(x)=LT_th(x),
\qquad
\forall h \in L^2({\mathbb R}^k,\Phi).
\label{eq1.8}
\ena
Now on the space $L^2({\mathbb R}^k,\Phi)$ (where $\Phi$ = $N(0,
{\mathbb I}_k )$ is the $k$-dimensional
standard Normal ), $L$ is self adjoint and has a spectral gap, with the eigenvalue $0$
corresponding to the invariant distribution $\Phi$ (or the constant function ${\mathbf 1}$
on $L^2({\mathbb R}^k,\Phi)$).
This may be deduced from the fact that the Normal density
$p(t;x,y)$
(with mean vector
$e^{-t}x$ and dispersion matrix $(1-e^{-2t}){\mathbb I}_k$)
converges to the standard Normal density $\phi(y)$ exponentially fast as
$t \rightarrow \infty$, for every initial state $x$. Else, one can compute
the set of eigenvalues of $L$, namely $\{0,-1,-2,\ldots \}-$ with eigenfunctions expressed in terms
of Hermite polynomials \cite[page 487]{bhattacharyawaymire}. In particular, $L^{-1}$ is
a bounded operator on
$1^{\bot}$  and is given by
\bea
L^{-1}\tilde{h}=
-\int_0^\infty
T_s \tilde{h}(x)ds, 
\qquad \forall \tilde{h} =h-\int h d\Phi\in L^2({\mathbb R}^k,\Phi).
\label{eq1.9}
\ena
To check this, note that by (\ref{eq1.8})
\bea
\tilde{h}=
-\int_0^\infty
\frac{\partial}{\partial s}
T_s \tilde{h}(x)ds
= 
-\int_0^\infty
LT_s\tilde{h}(x)ds
=L\left(
-\int_0^\infty
T_s\tilde{h}(x)ds
\right).
\ena
For our purposes $h={\mathbf 1}_C$: the indicator function of a Borel convex subset $C$
of ${\mathbb R}^k$.

A \underline{smooth approximation} of $\tilde{h}$ is $T_t\tilde{h}$
for small $t>0$ (since $T_t\tilde{h} $ is infinitely differentiable).
Also, by (\ref{eq1.9})
\bea
\psi_t(x)\equiv
L^{-1}T_t\tilde{h} (x)&=&
-\int_0^\infty
T_sT_t\tilde{h}(x)ds
=-\int_0^\infty
T_{s+t}\tilde{h}(x)ds
=
-\int_t^\infty
T_s\tilde{h}(x)ds
\label{eq1.10}
\\
&=&
-\int_t^\infty
\{
\int_{{\mathbb R}^k}
\tilde{h}
(e^{-s}x+\sqrt{1-e^{-2s}}z
)\phi(z)dz
\}
ds
\nonumber
\ena
where $\phi$ is the $k$-dimensional standard Normal  density.
We have expressed $T_s\tilde{h}(x)\equiv E[\tilde{h}(X_s)|X_0=x]$
in (\ref{eq1.10}) as
\bea
E[\tilde{h}(X_s)|X_0=x]=
E\tilde{h}
(e^{-s}x+\sqrt{1-e^{-2s}}Z),
\label{eq1.11}
\ena
where $Z$ is a standard Normal  $N(0,
{\mathbb I}_k )$. For $X_s$ has the same distribution
as $e^{-s}x+\sqrt{1-e^{-2s}}Z$. Now note that using
(\ref{eq1.10}), one may write
\bea
\hskip0.3cm
T_t\tilde{h}(x)=
L(L^{-1}T_t\tilde{h}(x))
=
\Delta
(L^{-1}T_t\tilde{h}(x))-
x \bullet
\nabla
(L^{-1}T_t\tilde{h}(x))=
\Delta \psi_t(x)-x\bullet
\nabla \psi_t(x).
\label{eq1.12}
\ena 

For the problem at hand
(see \ref{eq1.1})
 $Q_0=\Phi$ and $Q=Q_{(n)}$ is the distribution of $S_n=
\frac{1}{\sqrt{n}} (Y_1+Y_2+\cdots +Y_n)=
(X_1+X_2+\cdots +X_n)$, $(X_j=Y_j/\sqrt{n})$,
where $Y_j's$ are i.i.d. mean-zero with covariance
matrix ${\mathbb I}_k$
and finite absolute third moment
$$\rho_3=E||Y_1||^3=E(\sum_{i=1}^k(Y_1^{(i)})^2)^\frac{3}{2}.$$
We want to estimate 
\bea
E\tilde{h}(S_n)=Eh(S_n)-\int h d\Phi
\label{eq1.15}
\ena
for $h={\mathbf 1}_C$, $C \in {\mathcal C}-$the class
of all Borel convex sets in ${\mathbb R}^k$.

For this we first estimate (see  (\ref{eq1.12})), for small $t>0$,
\bea
ET_t\tilde{h}(S_n)=
E\left[\Delta \psi_t(S_n)-S_n\bullet \nabla \psi_t (S_n)\right]
\label{eq1.16}
\ena
This is done in Section 3.
The next step is to estimate, for small $t>0$,
\bea
ET_t\tilde{h}(S_n)-
E\tilde{h}(S_n)
\label{eq1.17}
\ena
which is carried out in Section 4.
Combining the estimates of (\ref{eq1.16}) and (\ref{eq1.17}),
and with a suitable choice of $t>0$, one arrives at the desired estimation of
(\ref{eq1.15}).

We will write
\bea
\delta_n=\sup_{\{h={\mathbf 1_C}:C\in{\mathcal C}\}} |\int h dQ_{(n)}-\int h d\Phi|.
\ena

\section{Derivatives of $\psi_t\equiv L^{-1}T_t \tilde{h}$}
Before we engage in the estimation of (\ref{eq1.16}) and (\ref{eq1.17}), it is useful
to compute certain derivatives of $\psi_t$.
$$\mbox{ Let }D_i=\frac{\partial}{\partial x_i},\qquad
D_{ii'}=\frac{\partial^2 }{\partial x_i \partial x_{i'}},
\qquad
D_{ii'i''}=\frac{\partial^3 }{\partial x_i \partial x_{i'}\partial x_{i''}}, \mbox{ etc
..}$$
Then, using (\ref{eq1.10}),
\bea
\hskip-1cm
D_i\psi_t(x)
&=&
-\int_t^\infty
\left[
\int_{{\mathbb R}^k}
\tilde{h} (y) (2\pi(1-e^{-2s}))^{-\frac{k}{2}}
\frac{(y_i-e^{-s}x_i)}{1-e^{-2s}}
e^{-s}
\cdot 
\exp\{- \frac{||y-e^{-s}x ||^2}{2(1-e^{-2s})}
\} dy \right] ds
\nonumber \\
&=&
-\int_t^\infty
\left[
\int_{{\mathbb R}^k}
\tilde{h} (y) (2\pi(1-e^{-2s}))^{-\frac{k}{2}}
\frac{e^{-s}}{\sqrt{1-e^{-2s}}}
\cdot \frac{(y_i-e^{-s}x_i)}{\sqrt{1-e^{-2s}}}
\cdot \exp\{ -\frac{||y-e^{-s}x ||^2}{2(1-e^{-2s})}
\}
dy
\right]
ds
\nonumber
\\
&=&
-\int_t^\infty
(\frac{e^{-s}}{\sqrt{1-e^{-2s}}})
\left[
\int_{{\mathbb R}^k}
\tilde{h} (e^{-s}x+\sqrt{1-e^{-2s}}z)z_i\phi(z)
dz
\right]
ds,
\qquad
z_i\phi(z)=-\frac{\partial}{\partial z_i}\phi(z))=-D_i\phi(z)
\ena
{using the change of variables}
$$
z=\frac{y-e^{-s}x}{\sqrt{1-e^{-2s}}}.
$$
In the same manner, one has, using $D_{i}$,$D_{ii'}$, etc for derivatives
$\frac{\partial}{\partial z_i}$,
$\frac{\partial^2}{\partial z_i\partial z_{i'}}$,
etc,
\bea
D_{ii'}\psi_t(x)
&=&
-\int_t^\infty
(\frac{e^{-s}}{\sqrt{1-e^{-2s}}})^2
\left[
\int_{{\mathbb R}^k}
\tilde{h} (e^{-s}x+\sqrt{1-e^{-2s}}z)
\cdot D_{ii'}
\phi(z)
dz
\right]
ds,
\nonumber\\
D_{ii'i''}\psi_t(x)
&=&
-\int_t^\infty
(\frac{e^{-s}}{\sqrt{1-e^{-2s}}})^3
\left[
\int_{{\mathbb R}^k}
\tilde{h} (e^{-s}x+\sqrt{1-e^{-2s}}z)
\cdot (-D_{ii'i''}
\phi(z))
dz
\right]
ds.
\ena
The following estimate is used in the next section:
\bea
\sup_{u\in{\mathbb R}^k}
\left|
\int_{{\mathbb R}^k}\int_{{\mathbb R}^k}
\tilde{h}(\sqrt{\frac{n-1}{n}}e^{-s}x+e^{-s}
u
+\sqrt{1-e^{-2s}}z)
\phi(x)
D_{ii'i''}\phi(z)dxdz \right|
\label{eq2.3}
\\
\nonumber
\leq c_0ke^{2s}(1-e^{-2s}).
\ena
To prove this, write 
$a=\sqrt{\frac{n}{n-1}}e^{s}\sqrt{1-e^{-2s}}$
and change variables
$x \longrightarrow y=x+az$. Then
\bea
\label{eq2.4}
\phi(x)=\phi(y-az)=\phi(y)-az\cdot\nabla\phi(y)+a^2
\sum_{r,r'=1}^kz_rz_{r'}
\int_0^1 (1-v)D_{rr'}\phi(y-vaz)dv,
\ena
so that
\bea
\tilde{h}(\sqrt{\frac{n-1}{n}}e^{-s}x+
e^{-s}u+
\sqrt{1-e^{-2s}}z)
=
\tilde{h}(\sqrt{\frac{n-1}{n}}e^{-s}y+
e^{-s}u),
\nonumber
\ena
and the double integral in (\ref{eq2.3}) becomes
\begin{align}
\label{eq2.5}
&\\
\int_{{\mathbb R}^k}\int_{{\mathbb R}^k}
\tilde{h}(\sqrt{\frac{n-1}{n}}e^{-s}y+e^{-s}
u)
\left[
\phi(y)-az\cdot\nabla\phi(y)+a^2
\sum_{r,r'=1}^kz_rz_{r'}
\int_0^1 (1-v)D_{rr'}\phi(y-vaz)dv
\right]
\bullet
D_{ii'i''}\phi(z)dzdy
\nonumber
\end{align}

Note that the integrals of $D_{ii'i''}\phi(z)$ and $z_{i_0}D_{ii'i''}\phi(z)$
vanish for $i,i',i'',i_0$, so that 
\bea
\int_{{\mathbb R}^k}
\tilde{h}(\sqrt{\frac{n-1}{n}}e^{-s}y+e^{-s}u)
(\phi(y)-az\cdot\nabla\phi(y))
D_{ii'i''}\phi(z)dz=0
\ena
The magnitude of the last term  on the right in (\ref{eq2.4}) is
\bea
&&
\left|
a^2
\int_0^1 (1-v)
\left[
\sum_{r,r'=1}^kz_rz_{r'}
(y-vaz)_r 
(y-vaz)_{r'}
-\sum_{r=1}^kz_r^2
\right]
\phi(y-avz)dv 
\right|\\
&&\leq
a^2
\int_0^1 (1-v)
\left[
\sum_{r,r'=1}^kz_rz_{r'}
(y-vaz)_r (y-vaz)_{r'}
+\sum_{r=1}^kz_r^2
\right]
\phi(y-avz)dv ,
\nonumber
\ena
since the sum $\sum_{r,r'}$ above is nonnegative. Bounding $|\tilde{h}|$
by 1,
it follows from (\ref{eq2.5})-(2.7)
that
 the left side of (\ref{eq2.3}) is  no more than
\begin{align}
&a^2\!
\int_0^1\! (1-v)\!
\left\{
\!
\int_{{\mathbb R}^k}
\left[
\sum_{r\neq r'}^kz_rz_{r'}
\int_{{\mathbb R}^k}
(y-vaz)_r (y-vaz)_{r'}
\phi(y-vaz)dy
\!
+
\!
\sum_{r=1}^kz_r^2
\int_{{\mathbb R}^k}
\{\!
(y-vaz)_r^2
+1\}\!
\phi(y-avz)dy
\right]
\!
|D_{ii'i''}\phi(z)|dz
\right\}dv \nonumber \\
=&
a^2
\int_0^1 (1-v)
\left\{
\int_{{\mathbb R}^k}
2
\sum_{r=1}^kz_r^2
|D_{ii'i''}\phi(z)|dz\!
\right\}
\!dv,
\label{eq2.8}
\end{align}
from which (\ref{eq2.3}) follows.
\section{Estimation of $T_t\tilde{h}(S_n)$}
By  (\ref{eq1.12}),
\bea
T_t\tilde{h}(S_n)=
L(L^{-1}T_t\tilde{h})(S_n)
=
L\psi_t(S_n)=
\Delta \psi_t(S_n)-S_n\bullet
\nabla
\psi_t(S_n)
\ena
Consider the Taylor expansions
$$\Delta\psi_t(S_n)\equiv \sum_{i=1}^k D_{ii}\psi_t(S_n)=
\sum_{i=1}^k D_{ii}\psi_t(S_n-X_1)
+\sum_{i,i'=1}^k \int_0^1
X_1^{(i')}
D_{iii'}
\psi_t(S_n-X_1+v X_1)
dv,$$
\begin{align}
S_n\cdot \nabla \psi_t(S_n)&=\sum_{j=1}^nX_j\cdot\nabla
\psi_t(S_n)=
\sum_{j=1}^n
\sum_{i=1}^k
X_j^{(i)}\cdot D_i \psi_t(S_n)
\nonumber
\\
&=
\sum_{j=1}^n
\left[
\sum_{i=1}^k
X_j^{(i)}\cdot D_i \psi_t(S_n-X_j)+
\sum_{i,i'=1}^k
X_j^{(i)}
X_j^{(i')}
D_{ii'} \psi_t(S_n-X_j)
+
\right.
\\
+& 
\left.
\sum_{i,i',i''=1}^k
X_j^{(i)}
X_j^{(i')}
X_j^{(i'')}
\int_0^1
(1-v )
D_{ii'i''} \psi_t(S_n-X_j+v X_j)
dv 
\right]
\nonumber
\end{align}

Recalling that $X_j=\frac{Y_j}{\sqrt{n}}$,
$EY_j=0$,
$EX_j^{(i)}X_j^{(i')}=
\frac{1}{n}EY_j^{(i)}Y_j^{(i')}
=\frac{1}{n}\delta_{ii'}$
and 
$X_j$ and $S_n-X_j$ are independent,
\begin{align}
E\Delta\psi_t(S_n)
=
E\left[
\sum_{i=1}^k
D_{ii} \psi_t(S_n-X_1)\right]+
E\left[
\sum_{i,i'=1}^k
\frac{Y_1^{(i)}}{\sqrt{n}}
\int_0^1
D_{iii'} \psi_t(S_n-X_1+v X_1)
dv 
\right],
\end{align}

\begin{align}
ES_n\cdot \nabla \psi_t(S_n)
=
E\left[
\sum_{i=1}^k
D_{ii} \psi_t(S_n-X_1)\right]+
\frac{1}{\sqrt{n}}
\sum_{i,i',i''=1}^k
E\left[
Y_1^{(i)}
Y_1^{(i')}
Y_1^{(i'')}
\int_0^1
(1-v )
D_{ii'i''} \psi_t(S_n-X_1+v X_1)
dv 
\right].
\end{align}
Hence
\begin{align}
ET_t\tilde{h}(S_n)
=
E\left[
\sum_{i,i'=1}^k
\frac{Y_1^{(i)}}{\sqrt{n}}
\int_0^1
D_{iii'} \psi_t(S_n-X_1+v X_1)
dv 
-\frac{1}{\sqrt{n}}
\sum_{i,i',i''=1}^k
Y_1^{(i)}
Y_1^{(i')}
Y_1^{(i'')}
\int_0^1
(1-v )
D_{ii'i''} \psi_t(S_n-X_1+v X_1)
dv 
\right].
\label{eq3.4}
\end{align}
One may then write
\bea
ET_t\tilde{h}(S_n)=E
[E(\bullet \bullet \bullet|Y_1)] 
\label{eq3.5}
\ena
where $\bullet \bullet \bullet$ is the quantity
within square brackets in (\ref{eq3.4}),
i.e.,
\begin{align}
E[T_t\tilde{h}(S_n)|Y_1]&=
\frac{1}{\sqrt{n}}
\sum_{i,i'=1}^k
Y_1^{(i')}
\int_0^1
E\left[
D_{iii'} \psi_t(S_n-X_1+v X_1)
|Y_1
\right]
dv
\label{eq3.6}
\\
&
-\frac{1}{\sqrt{n}}\sum_{i,i',i''=1}^k
Y_1^{(i)}
Y_1^{(i')}
Y_1^{(i'')}
\int_0^1(1-v)E[D_{ii'i''}\psi_t(S_n-X_1+vX_1)|Y_1]dv
\nonumber
\end{align}
The first term on the right side in (\ref{eq3.6}) equals
\begin{align}
&\frac{1}{\sqrt{n}}
\sum_{i,i'=1}^k
Y_1^{(i')}
\left(
-\int_t^\infty
(\frac{e^{-s}}{\sqrt{1-e^{-2s}}})^3
\int_0^1
\left\{
E\left[
\int_{{\mathbb R}^k}
\tilde{h} (e^{-s}(S_n-X_1)+e^{-s}v X_1+
\sqrt{1-e^{-2s}}z)
\cdot (-D_{iii'}
\phi(z))
dz
\vert
Y_1
\right]
\right\}
dv 
\right)
ds
\nonumber
\\
&=
\frac{1}{\sqrt{n}}\!
\sum_{i,i'=1}^k
Y_1^{(i')}
\!
\int_t^\infty
(\frac{e^{-s}}{\sqrt{1-e^{-2s}}})^3
\left(
\int_0^1
\left\{
\int_{{\mathbb R}^k}
\left[
\int_{{\mathbb R}^k}
\tilde{h} (e^{-s}
\sqrt{\frac{n-1}{n}}
x+e^{-s}v X_1+
\sqrt{1-e^{-2s}}z)
dQ_{(n-1)}(x)
\right]\!
D_{iii'}
\phi(z)\!
dz\!
\right\}
dv \!
\right)\!
ds,
\label{eq3.8}
\end{align}
noting that the distribution of $S_n-X_1=\sqrt{\frac{n-1}{n}}(\frac{Y_2+Y_3+\cdots
+Y_n}{\sqrt{n-1}})$ is that of
$\sqrt{\frac{n-1}{n}}V$, where $V$ has distribution $Q_{(n-1)}$. Therefore,
(\ref{eq3.8}) is equal to
\begin{align}
&
\frac{1}{\sqrt{n}}
\sum_{i,i'=1}^k
Y_1^{(i')}
\int_t^\infty
(\frac{e^{-s}}{\sqrt{1-e^{-2s}}})^3
\int_0^1
\left\{
\int_{{\mathbb R}^k}
\left[
\int_{{\mathbb R}^k}
\right.
\right.
&
\left.
\tilde{h} (e^{-s}
\sqrt{\frac{n-1}{n}}
x+e^{-s}v X_1+
\sqrt{1-e^{-2s}}z)
(d(Q_{(n-1)}(x)-\Phi(x))+d\Phi(x))
\right]
\label{eq3.9}
\\
&&
\cdot
D_{iii'}
\phi(z)
dz
\bigg\}
dv 
ds
\nonumber
\end{align}

Since the class of functions $h={\mathbf 1}_C$, where $C$ ranges over all
Borel convex subsets of ${\mathbb R}^k$, is invariant under translation,
and $bC$ is convex if $C$ is convex
($bC=\{ bx:x \in C\}, \forall b >0$),
\bea
\left|
\int_{{\mathbb R}^k}
\tilde{h} (e^{-s}
\sqrt{\frac{n-1}{n}}
x+e^{-s}v X_1+
\sqrt{1-e^{-2s}}z)
(d(Q_{(n-1)}(x)-\Phi(x)))
\right|
\leq
\delta_{n-1}.
\label{eq3.10}
\ena
Similarly, the second term on the right 
in (\ref{eq3.6}) equals
\begin{align}
&-
\frac{1}{\sqrt{n}}
\sum_{i,i',i''=1}^k
Y_1^{(i)}
Y_1^{(i')}
Y_1^{(i'')}
\int_t^\infty
(\frac{e^{-s}}{\sqrt{1-e^{-2s}}})^3
\int_0^1
\left\{
\int_{{\mathbb R}^k}
\left[
\int_{{\mathbb R}^k}
\right.
\right.
\tilde{h} (e^{-s}
\sqrt{\frac{n-1}{n}}
x+e^{-s}v X_1+
\sqrt{1-e^{-2s}}z)
\label{eq3.11}
\\
&
\bigg(
d(Q_{(n-1)}(x)-\Phi(x))+d\Phi(x)
\bigg)
\bigg]
. D_{ii'i''} 
\phi(z)dz 
\bigg\}
dv 
\nonumber
\end{align}
Again, the inner integral in (\ref{eq3.11})
with regard to $Q_{(n-1)}-\Phi$ is estimated
by (\ref{eq3.10}). Therefore,
using (\ref{eq2.3}) for the remaining integration
with regard to $\Phi$ in (\ref{eq3.8}), (\ref{eq3.11}). 
\begin{align}
\left| ET_t\tilde{h}(S_n) \right|
&\leq
\frac{1}{\sqrt{n}}
\sum_{i,i'=1}^k
E|Y_1^{(i')}|
\left(
\int_t^\infty
(\frac{e^{-s}}{\sqrt{1-e^{-2s}}})^3
\left[\delta_{n-1}
\int_{{\mathbb R}^k}
|D_{iii'}\phi(z)|dz
+c_0ke^{2s}(1-e^{-2s})
\right]
ds
\right)
\nonumber
\\
&\!+\!
\frac{1}{\sqrt{n}}\!
\sum_{i,i',i''=1}^k
\!
E\left|Y_1^{(i)}
Y_1^{(i')}
Y_1^{(i'')} \right|
\left(
\int_t^\infty
(\frac{e^{-s}}{\sqrt{1-e^{-2s}}})^3
\left[\delta_{n-1}
\left(\int_0^1
(1-v )dv\right)
\cdot
\int_{{\mathbb R}^k}
|D_{ii'i''}\phi(z)|dz\!
+c_0ke^{2s}(1-e^{-2s})
\right]\!
ds \! \right).
\label{eq3.12}
\end{align}
Next, the first two terms on the right in (\ref{eq3.12}) may be estimated
by using
\bea
\int_{{\mathbb R}^k}
|D_{ii'i''} \phi(z)| dz=
\begin{cases}
E|(Z_1^{(i)})^2-1|\cdot E|Z_1^{(i')}| \leq 1
&\forall i \neq i',i''=i or i',\\
E|(Z_1^{(i)})^3-Z^{(i)}| \leq \sqrt{6}
&\forall i = i'=i'',\\
\end{cases}
\nonumber
\\
\int_{{\mathbb R}^k}
|D_{ii'i''} \phi(z)| dz=
E|Z_1^{(i)}Z_1^{(i')}Z_1^{(i'')}
|\leq 1
\mbox{ if }
i,i',i''
\mbox{ are all distinct.}
\label{eq3.13}
\ena
Finally, note that  
\begin{equation}
\frac{e^{-s}}{\sqrt{1-e^{-2s}}}
\leq \frac{1}{\sqrt{2s}} \qquad (s>0), 
\label{eq3.14}
\end{equation}
 so that
\begin{align}
\int_0^\infty
\frac{e^{-s}}{\sqrt{1-e^{-2s}}} ds
=c_0' < \infty,
\qquad
\int_t^\infty
(\frac{e^{-s}}{\sqrt{1-e^{-2s}}})^3 ds
\leq (2t)^{-\frac{1}{2}}.
\label{eq3.15}
\end{align}
Hence, using (\ref{eq3.12})-(\ref{eq3.15}),
together with the estimates 
$$
E\sum_{i,i',i''=1}^k
\left|
Y_1^{(i)}
Y_1^{(i')}
Y_1^{(i'')}
\right|
\leq
k^{\frac{3}{2}} \rho_3, \quad 
E\sum_{i,i'=1}^k|Y_1^{(i')}|
\leq k^{\frac{1}{2}}\rho_3,
$$
one has
\bea
|ET_t\tilde{h}(S_n)|\leq
c_1k^{3/2}\rho_3(
\frac{\delta_{n-1}}{\sqrt{n}\sqrt{t}})
+
\frac{c_2k^{5/2}\rho_3}{n^{1/2}} 
\label{eq3.16}
\ena
\section{The smoothing inequality and the Estimation of $\delta_n$}
Let ${\mathcal H}=\left\{ {\mathbf 1}_C, C \in {\mathcal C}
\right\}$, where $ {\mathcal C}$ is the class of all Borel convex subsets of
${\mathbb R}^k$.
As before, $\tilde{h}=h-\int h d\Phi$. We also write $G_b$ as the distribution of $bW$, if $W$ has distribution $G (b>0)$. Recall that
(see \ref{eq1.11})
$T_t\tilde{h}(x)=E\tilde{h}
(e^{-t}x+\sqrt{1-e^{-2t}}Z)$, where $Z$ has the standard Normal
distribution $\Phi=N(0,{\mathbb I}_k)$, which we take to be independent of $S_n$.
Then
\begin{align}
ET_t\tilde{h}(S_n)
&=
E
\tilde{h}
(e^{-t}S_n+\sqrt{1-e^{-2t}}Z
)
=
\int_{{\mathbb R}^k}
\int_{{\mathbb R}^k}
\tilde{h}
(e^{-t}x+\sqrt{1-e^{-2t}}z
)
dQ_{(n)}(x)
\phi(z)dz
\nonumber
\\
&=
\int_{{\mathbb R}^k}
\tilde{h}
d((Q_{(n)})_{e^{-t}}
\star
\Phi_{\sqrt{1-e^{-2t}}}
)
=
\int_{{\mathbb R}^k}
\tilde{h}
d((Q_{(n)})_{e^{-t}}
-
\Phi_{e^{-t}})
\star
\Phi_{\sqrt{1-e^{-2t}}}
\label{eq4.2}
\end{align}
The introduction of the extra term
$\Phi_{e^{-t}}
\star
\Phi_{\sqrt{1-e^{-2t}}}=\Phi$
does not affect the integration in the last step
since
$\int_{{\mathbb R}^k}
\tilde{h}
d\Phi=0.
$

Since the last integration is with respect to the difference
between two probability 
measures, its value is unchanged if we replace
$\tilde{h}$ by $h$.
Hence
\begin{align}
ET_t\tilde{h}(S_n)
=\int_{{\mathbb R}^k}
h
d[Q_{(n)})_{e^{-t}}
-
\Phi_{e^{-t}}]
\star
\Phi_{\sqrt{1-e^{-2t}}}
\; .
\label{eq4.3}
\end{align}
Also the class ${\mathcal C}$ is invariant under multiplication 
$C\longrightarrow bC$ where $b>0$ is given. Therefore,
\bea
\delta_n=\sup_{h \in {\mathcal H}} | 
E\tilde{h}(S_n)|=
\sup_{h \in {\mathcal H}} | 
\int h d(Q_{(n)}-\Phi)|
=
\sup_{h \in {\mathcal H}} | 
\int h d\left[
(Q_{(n)})_{e^{-t}}
-
\Phi_{e^{-t}}
\right].
\label{eq4.4}
\ena
Thus (\ref{eq4.3}) is a perturbation (or, smoothing) of the integral
in (\ref{eq4.4}) by convolution
with $\Phi_{\sqrt{1-e^{-2t}}}$. If $\epsilon>0$
is a constant such that
\bea
\Phi_{\sqrt{1-e^{-2t}}}
\left(\left\{
|z|< \epsilon
\right\}
\right)=\frac{7}{8},
\label{eq4.5}
\ena
then the smoothing inequality below applies,
with $\mu=(Q_{(n)})_{e^{-t}}$,
 $\nu =\Phi_{e^{-t}}, K=\Phi_{\sqrt{1-e^{-2t}}}, f=h={\mathbf 1}_C$, $\alpha=7/8$,
 and $\epsilon$ as in (\ref{eq4.5}).

\subsection*{Smoothing Inequality}
Let $\mu,\nu , K$ be probability measures on
${\mathbb R}^k$,
$K(\left\{x:
|x|< \epsilon
\right\}
)=\alpha
>\frac{1}{2}
$.Then for every bounded measurable $f$ one has
\bea
\left|
\int_{{\mathbb R}^k}
f d(\mu- \nu ) \right|
\leq (2 \alpha -1)^{-1}
\left[\gamma^*(f:\epsilon)+\omega_f^*(2\epsilon:\nu )
\right]
\label{eq4.6}
\ena
where, 
\beas
f^+_\epsilon(x)&=&
\sup \{ f(y): |y-x|<\epsilon \},
f^-_\epsilon(x)=
\inf \{ f(y): |y-x|<\epsilon \},\\
\mbox{
and }
\gamma(f:\epsilon)&=&\max\left\{
\int_{{\mathbb R}^k}
\left|f^+_\epsilon d(\mu- \nu) \right|
,
\int_{{\mathbb R}^k}
\left|f^-_\epsilon d(\mu- \nu ) \right|
\right\},\\
\gamma^*(f:\epsilon)&=&\sup_{y\in
{\mathbb R}^k}
\gamma(f_y:\epsilon),\qquad 
f_y(x)\equiv
f(x+y),
\\
\omega_f(x:\epsilon)&=&
\sup \{ |f(y)-f(x)|: |y-x|<\epsilon \},
\quad
\omega_f(\epsilon:v )=
\int
\omega_f(x:\epsilon)d\nu (x),\\
\omega^*(f:\epsilon)&=&\sup_{y\in
{\mathbb R}^k}
\omega_{f_y}(\epsilon:\nu ).
\enas

For a proof of the inequality (\ref{eq4.6}) see
Bhattacharya and Rao
\cite{bhattacharyarao}, Lemma 11.4.
With $h={\mathbf 1}_C$ one gets 
$h^+_{\epsilon}={\mathbf 1}_{C^\epsilon}$
$h^-_{\epsilon}={\mathbf 1}_{C^{-\epsilon}}$, where
$
{C^\epsilon}
=\left\{x:
dist(x,C)< \epsilon
\right\}
$
,${C^{-\epsilon}}=\left\{x: \mbox{open ball of radius }
\epsilon \mbox{ and center }
x \mbox{ is contained in } C\right\}$
are both convex, so that 
$$
\gamma(h:\epsilon)\leq \max \left\{
\int
{\mathbf 1}_{C^\epsilon}
d 
\left[
(Q_{(n)})_{e^{-t}}
-\Phi_{e^{-t}}
\right]
\star \Phi_{\sqrt{1-e^{-2t}}},\;
\int {\mathbf 1}_{C^{-\epsilon}}
d \left[
(Q_{(n)})_{e^{-t}}
-\Phi_{e^{-t}}
\right]
\star \Phi_{\sqrt{1-e^{-2t}}}
\right\}
\leq \sup_{h\in {\mathcal H}} \mid ET_t\tilde{h}(S_n)\mid
.
$$
Since ${\mathcal C}$ is invariant under translation one then obtains
\bea
\gamma^*(h:\epsilon)\leq
\sup_{h\in {\mathcal H}} \mid ET_t\tilde{h}(S_n)\mid.
\label{eq4.7}
\ena
Also, letting $Z$ be standard Normal $N(0,\one_k)$,
\bea
\omega^*_h(2\epsilon:\Phi_{e^{-t}})&=&P(e^{-t}Z \in (\partial C)^{2\epsilon})
\nonumber\\
&=&P(Z \in e^{t}(\partial C)^{2\epsilon})
\leq c_3\sqrt{k}2\epsilon e^{t}.
\label{eq4.8}
\ena
where $c_3>0$ is a constant (see Bhattacharya and Rao \cite{bhattacharyarao},
Theorem 3.1).
From (\ref{eq4.5})
one gets 
$$P\left(\left|
\sqrt{1-e^{-2t}}
Z\right|< \epsilon
\right)=\frac{7}{8},
\qquad
P\left(\left|
Z\right|< \frac{\epsilon}{\sqrt{1-e^{-2t}}}
\right)=\frac{7}{8}
$$
so that $\epsilon/\sqrt{1-e^{-2t}}=a_k$,
where $a_k$ satisfies 
$P(|Z|<a_k)=\frac{7}{8}$.
It is simple to check that $a_k=O(\sqrt{k})$, as $k\longrightarrow \infty$,
and 
\bea
a_k\leq c_4\sqrt{k}, \epsilon=a_k\sqrt{1-e^{-2t}}\leq
c_4\sqrt{k}\sqrt{1-e^{-2t}}\leq
 c_4\sqrt{k}\sqrt{2t}
\ena
Using this estimate of $\epsilon$ in (\ref{eq4.8}), one obtains
\bea
\omega^*_h(2\epsilon:\Phi_{e^{-t}})\leq
c_5k\sqrt{t}e^t
\label{eq4.8}
\ena
The smoothing inequality now yields (use (\ref{eq4.4}),(\ref{eq4.7}),
(\ref{eq4.8}) in (\ref{eq4.6}))
\bea
\delta_n \leq 
\frac{4}{3}
\left[
\sup_{h\in {\mathcal H}} \mid ET_t\tilde{h}(S_n)\mid
+ c_5k\sqrt{t}e^t
\right]
\label{eq4.10}
\ena
Now use (\ref{eq3.14}) in (\ref{eq4.11}) to get
\bea
\delta_n \leq
(c_6k^{3/2}\rho_3)\frac{\delta_{n-1}}{\sqrt{n}\sqrt{t}}+
\frac{c_7 k^{5/2} \rho_3}{n^{1/2}}+c_8k\sqrt{t}e^t.
\label{eq4.11}
\ena
By comparing 
the first and third terms on the right, an optimal order of $t$
is obtained as 
$$
t=\min \left\{1,\frac{\sqrt{k}\delta_{n-1}\rho_3}{\sqrt{n}} 
\right\}.
$$
It follows that
\bea
\delta_n \leq
(c_9k^{5/4}\rho_3^{1/2})\frac{\delta_{n-1}^{1/2}}{n^{\frac{1}{4}}}+
\frac{c_7 k^{3/2} \rho_3}{n^{1/2}}.
\label{eq4.12}
\ena
Consider now the induction hypothesis : The inequality
\bea
\delta_n\leq \frac{ck^{5/2}}{\sqrt{n}} \rho_3
\label{eq4.13}
\ena
holds for some $n\geq 1$ and an absolute constant $c\geq 1$ specified below.
Note that (\ref{eq4.13}) clearly holds for $n\leq c^2k^5   \rho_3^2$
Since $c^2k^5   \rho_3^2 >k^8$,
suppose then (\ref{eq4.13})
holds for some $n=n_0\geq k^8$.
Then by (\ref{eq4.12})
 we can take $n_0\geq k^3$: under the induction
hypothesis, and (\ref{eq4.12}),
\begin{align}
\delta_{n_0+1} &\leq 
\frac{c_9\sqrt{c}k^{\frac{5}{4}+\frac{5}{4}\rho_3}}{((n_0(n_0+1))^{\frac{1}{4}}} 
+ \frac{c_7k^{3/2}\rho_3}{(n_0+1)^{\frac{1}{2}}}
\nonumber
\\
&\leq
\frac{c_{10}\sqrt{c}k^{\frac{5}{2}} \rho_3}{(n_0+1)^{\frac{1}{2}}} 
+ \frac{c_{7}k^{5/2}\rho_3}{2^9(n_0+1)^{\frac{1}{2}}}
\qquad
(c_{10}=c_9+1,
\frac{k^{-1}}{n_0+1}\leq k^{-9}\leq 2^{-9}, \quad \mbox{ for k }\geq 2).
\nonumber
\\
\end{align}
Now, choose $c$ to be the greater  of 1 and the  positive solution of  $c=c_{10}\sqrt{c} +
c_7$, to check that (\ref{eq4.13}) holds for $n=n_0+1$.
Hence (\ref{eq4.13}) holds for all n.

We have proved the following result.
\begin{theorem}
There exists an absolute 
constant $c>0$ such that
\bea
\delta_n \leq
\frac{c k^{\frac{5}{2}}\rho_3}{\sqrt{n}}
\ena
\end{theorem}

\section{The Non-Identically Distributed Case}
For the general case considered in \cite{gotze}, $X_j$'s $(1\leq j \leq n)$
are independent with zero means and
$\sum_{j=1}^n Cov X_j={\mathbb I}_k$. Assume
\begin{equation}
\beta_3\equiv \sum_{1 \leq j\leq n} E||X_j||^3 < \infty
\end{equation} 
Let 
$\left\{
\bar{X}_j: 1 \leq j\leq n
\right\}$ be an independent copy of
$\left\{
X_j: 1 \leq j\leq n
\right\}$. Then, writing $S_n=\sum_{j=1}^n X_j$, as before,
\bea
E \sum_{i=1}^k D_{ii}\psi_t(S_n)&=&
E \sum_{j=1}^n \sum_{i,i'=1}^k D_{ii'}\psi_t(S_n)
\bar{X}_j^{(i)}\bar{X}_j^{(i')}
\nonumber
\\
&=&E\left[
\sum_{j=1}^n \sum_{i,i'=1}^k D_{ii'}\psi_t(S_n-X_j)
\bar{X}_j^{(i)}\bar{X}_j^{(i')}
+
\sum_{j=1}^n \sum_{i,i',i''=1}^k \bar{X}^{(i)}_j
\bar{X}^{(i')}_j
\bar{X}^{(i'')}_j
\int_0^1
D_{ii'i''}\psi_t(S_n-X_j+vX_j)dv
\right],
\label{eq5.2}
\ena
 and 
\begin{align}
&E
\left[S_n\cdot \nabla \psi_t(S_n)\right]
=
E\left[
\sum_{j=1}^n
X_j\cdot \nabla \psi_t(S_n)\right]
\label{eq5.3}
\\
\nonumber
&=
E\left[
\sum_{j=1}^n
\left\{
X_j\cdot \nabla \psi_t(S_n-X_j)
\!+\!
\sum_{i,i'=1}^k
X_j^{(i)}
X_j^{(i')}
D_{ii'}\psi_t(S_n-X_j) \!
+ \!
\!
\sum_{i,i',i''=1}^k
X_j^{(i)}
X_j^{(i')}
X_j^{(i'')}
\!
\int_0^1(1-v)
D_{ii'i''}\psi_t(S_n-X_j+vX_j)dv
\right\}
\!
\right]
\end{align}
Substracting (\ref{eq5.3}) from (\ref{eq5.2})
and noting that 
$$ 
 EX_j\cdot \nabla\psi_t(S_n-X_j)=0,$$ one obtains
\begin{align}
 &ET_t\tilde{h}(S_n)=
 E\left[ 
  \sum_{j=1}^n\sum_{i,i',i''=1}^k
  \bar{X}_j^{(i)}   \bar{X}_j^{(i')} \bar{X}_j^{(i'')}
    \int_0^1 D_{ii'i''}\psi_t(S_n-X_j+vX_j) dv
    \right.
      \nonumber \\
   & -
      \sum_{j=1}^n
      \sum_{i,i',i''=1}^k
  \bar{X}_j^{(i)}   \bar{X}_j^{(i')} \bar{X}_j^{(i'')}
    \int_0^1 
    (1-v)
    D_{ii'i''}\psi_t(S_n-X_j+vX_j) dv
\bigg]
\label{eq5.4}
  \end{align}
 The estimation of the conditional expectation of the integrals $\int_0^1$ in (\ref{eq5.4}),
 given $X_j$, proceeds as in Section 3 (with $X_j$ in place of $X_1$).
 The only significant change is in the normalization in the argument of $\tilde{h}$
 (see (\ref{eq3.8}) - (\ref{eq3.11})) where, writing
 $N_j$ as the positive square root
 of the inverse of $Cov(S_n-X_j)$,
 \begin{align}
 \label{eq5.5}
& E\left[\tilde{h}(e^{-s}(S_n-X_j)+e^{-s}vX_j+\sqrt{1-e^{-2s}}z|X_j\right]\\
 &=E\left[\tilde{h}(e^{-s}N_j^{-1}(N_j(S_n-X_j))+e^{-s}vX_j+\sqrt{1-e^{-2s}}z|X_j\right]
 \nonumber\\
&=
 \int_{{\mathbb R}^k}
\tilde{h}
(e^{-s}N_j^{-1}x+e^{-s}vX_j+
\sqrt{1-e^{-2s}}z)dQ_{(n-1),j}(x)
&\! \! \! \! \!\!=
 \int_{{\mathbb R}^k}
\tilde{h}
(e^{-s}N_j^{-1}(x+
N_je^s\sqrt{1-e^{-2s}}z)+e^{-s}vX_j
)dQ_{(n-1),j}(x),
\nonumber
 \end{align}
 where $Q_{(n)}$ denotes the distribution of
 $S_n=\sum_1^nX_j$, and $Q_{(n-1),j}$
 that of $N_j(S_n-X_j)$, which has mean zero, covariance
 ${\mathbb I}_k$. As in Section 3, the last integration is divided into two parts:
 $d(Q_{(n-1),j}-\Phi)(x)+d\Phi(x)$.
 Since the class of Borel convex sets is invariant under non-singular affine linear
 transformations, the integral with regards to 
 $Q_{(n-1),j}-\Phi$
is bounded by $\delta_{n-1}$. For the integral with regards to $\Phi$, we change variables
$x\longrightarrow y=x+A_jz$, where $A_j=e^{-s}\sqrt{1-e^{-2s}}N_j$. The estimation of the integral
now proceeds as in (\ref{eq2.3})$-$(\ref{eq2.8}), with scalar $a$ replaced by the matrix $A_j$.
The effect of this is simply to change the sum
$a^2\sum_{r,r'}z_rz_{r'}D_{rr'}\phi(y-vaz)
$ 
in (2.4)
to 
 $$
 \sum_{r,r'=1}^k (A_jz)_r
 (A_jz)_{r'}D_{rr'}\phi(y-vA_jz)
 $$
 Arguing as in (\ref{eq2.3})$-$(\ref{eq2.8})
 one arrives at the upper bound for (\ref{eq5.5}) given by
 $$
c_0''k \| A_j\|^2=c_0''ke^{2s}(1-e^{-2s})\|N_j\|^2\leq c_0''ke^{2s}(1-e^{-2s})(1-\beta_3^{\frac{2}{3}})^{-1},
 $$
 using
 \begin{align}
 \label{eq5.6}
& \|N_j\|^2=\|({\mathbb I_k-Cov X_j})^{-\frac{1}{2}} \|^2=\| {\mathbb I}_k-Cov X_j\|^{-1},\\
& \| {\mathbb I}_k-Cov X_j\|=\sup_{|u|=1} u\cdot ({\mathbb I}_k-Cov X_j) u
=\sup_{|u|=1} (1-E(u.X_j)^2)
\nonumber
\\
&\geq 1-E|X_j|^2\geq 1-(E|X_j|^3)^{\frac{2}{3}}\geq 1-\beta_3^{\frac{2}{3}}
\nonumber
 \end{align}
 and assuming 
 \begin{equation}
 \label{eq5.7}
 \beta_3<1
 \end{equation}
 Proceeding as in Section 4 one arrives at the bound:
 \begin{equation}
 \label{eq5.8}
\delta_n \leq ck^{\frac{5}{2}}\beta_3.
 \end{equation}
If one takes the absolute constant $c>1$, then the $\beta_3$ may be assumed
to be smaller or equal to $c^{-1}k^{-\frac{5}{2}}$, and
$
(1-\beta_3^{\frac{2}{3}})^{-1}\leq (1-\frac{1}{c^{\frac{2}{3}}})^{-1}=c'.
$
The induction argument is similar.\\
Remark: If one defines
\bea
\label{eq5.9}
\gamma_3\equiv \sum_{j=1}^n E(\sum_{i=1}^k|X_j^{(i)} |)^3,
\ena
then
\bea
\nonumber
\sum_{j=1}^n \sum_{i,i',i''=1}^k
E|X_j^{(i)}X_j^{(i')}X_j^{(i'')}|
= \gamma_3,
\ena
Since $\gamma_3$ now replaces $k^{\frac{3}{2}}\beta_3$
in the computations, it follows that
\bea
\label{eq5.10}
\delta_n\leq ck\gamma_3
\ena
Since, $\gamma_3\leq k^{\frac{3}{2}}\beta_3$,
(\ref{eq5.10}) provides a better bound than (\ref{eq5.8}) or (\ref{eq4.13}).

\end{document}